\documentclass[12pt]{article}

\usepackage{amssymb}
\usepackage{amsmath}

\usepackage{url}

\usepackage{color}
\usepackage{hyperref}
\usepackage{forest}
\usetikzlibrary{positioning}

\def\R{\mathcal{R}}

\def\D{\mathcal{D}}

\def\<{\langle}
\def\>{\rangle}

\def\ie{{\it i.e.}, }
\def\card{{\rm \#}}

\newcommand{\Set}{\textsc{Set}}
\newcommand{\Seq}{\textsc{Seq}}
\newcommand{\Cyc}{\textsc{Cyc}}

\usepackage{xspace}

\def\ashuff#1#2#3{
\kern 1pt \vrule height#1 \overline{\vrule height#3 width 0pt
\hskip#2} \rule{.3pt}{#1}\overline{\vrule height#3 width 0pt
\hskip#2} \rule{.3pt}{#1} \kern 1pt }

\def\binom#1#2{\left(#1\atop#2\right)}

\def\Carre3#1{\left[\begin{array}{ccc}#1\end{array}\right]}

\newtheorem{definition}{Definition} % amsthm only
\newtheorem{theorem}{Theorem}
\newtheorem{corollary}{Corollary}
\newtheorem{proposition}{Proposition}
\newtheorem{example}{Example}

\def\eg{\emph{e.g. }}
\def\etc{\emph{etc.}}
\def\labels{\mathtt{labels}}

\begin{document}
%\maketitle
\title{On recursively defined combinatorial classes and labelled trees}

\author{Ali Chouria\footnote{ali.chouria1@univ-rouen.fr. Laboratoire LITIS - EA 4108, Universit{\'e} de Rouen-Normandie , Avenue de l'Universit{\'e} - BP 8, 76801 Saint-\'Etienne-du-Rouvray, ISLAIB, University of Jendouba, Tunisia },
 Vlad-Florin Dr\u{a}goi\footnote{vlad-florin.dragoi@univ-rouen.fr. Department of Mathematics and Computer Sciences ''Aurel Vlaicu'' University of Arad, Romania }  and Jean-Gabriel Luque\footnote{jean-gabriel.luque@univ-rouen.fr. Laboratoire LITIS - EA 4108, Universit{\'e} de Rouen-Normandie, Avenue de l'Universit{\'e} - BP 8, 76801 Saint-\'Etienne-du-Rouvray Cedex}}
\maketitle

\begin{abstract}
%\selectlanguage{english} 
We define and prove isomorphisms between three combinatorial classes involving labeled trees.  We also give an alternative proof by means of generating functions.
\vskip 0.5\baselineskip

% Text of abstract in French
\end{abstract}

\noindent{\footnotesize {\bf Keywords:} {Combinatorial classes, Catalan numbers,  labeled trees, rooted trees, generating functions.}

% main text

\section{Introduction \label{Sec1}}

Based on the theory of combinatorial species (see \eg\cite{BLL}), Flajolet and Sedjewick \cite{Flajolet} wrote a reference book on combinatorial analysis. 
In particular in the first part, they provided a list of basic constructions for exponential generating functions. Mainly,  complex combinatorial structures are obtained by combining the following three combinatorial classes: $\Set, \Seq$, and $\Cyc.$ For instance, they described \emph{surjections} ($\Seq(\Set)$), \emph{set partitions} ($\Set(\Set)$), \emph{alignments} ($\Seq(\Cyc)$) and \emph{permutations} ($\Set(\Cyc)$).
In this article, we focus on the later one because it has the remarkable property of having the same generating function as the combinatorial class of sequences. More precisely, our starting point consists in giving an explicit bijection between the class of set-of-cycles and the class of sequences (see Section \ref{sec:setcyc-seq}). Our goal is to study an example of a class defined inductively by a combinatorial class equation. We chose the equation $\Set(\Cyc(\mathcal C))=\mathcal C$ because the underlying combinatorics reveal a world rich in interpretation and provide fruitful perspectives. In particular, this equation reveals an isomorphism between sets of necklaces of planar labelled trees, forests of labelled trees, and rooted labelled trees (see Section \ref{sec:trees}). 
At the heart of the equation we study here are the Catalan numbers. They are involved in the enumeration of numerous classes of combinatorial objects of prime importance in computer science, \eg  Dyck paths, binary trees, non-crossing partitions \etc \cite{Catalan,Novelli,Flajolet}. Notice that more than 60 possible enumerations are listed in \cite{Catalan} and the enumerated objects lead to several applications in computer science, such as sorting techniques  based on binary trees \cite{Knuth}. In our case, Catalan numbers are representing ordered trees. The article ends with Section \ref{sec:rec-tree-classes}, where we propose other recursive tree-like combinatorial classes.      

\section{Background and notations}
We recall here well known definitions and results concerning combinatorial classes and generating functions. The material contained in this section mainly refers to \cite{Flajolet}.
\subsection{Combinatorial classes}
%\paragraph{Definitions and properties}
In the most general context, a \emph{combinatorial class} is a triplet $(\mathcal O,\mathcal P,\omega)$ where $\mathcal O$ is the discrete set of the combinatorial objects we want to enumerate, $\mathcal P$ is the discrete set of the properties  in regard to which you want to enumerate our objects, and $\omega:\mathcal O\rightarrow \mathcal P$ is a map such that for every $p\in\mathcal P$ the preimage $\omega^{-1}(p)$ is finite, which is a minimal requirement in order to be able to enumerate the objects of $\mathcal O$ with respect to the properties of $\mathcal P$.\\
We consider the restricted context where $\mathcal P=\mathbb N$ and the preimage of $0$ contains only one element.
More formally, a combinatorial class is a pair  $\mathcal{C}=(\mathcal O,\omega)$ where $\omega : \mathcal{O} \rightarrow \mathbb{N}$ is such that ${\rm card}(\omega^{-1}(n))<\infty$ for any integer $n$. For the sake of simplicity, and when there is no ambiguity, we use the same name for a class and the set of its objects.
 Then, we denote by $|\mu|$ the \emph{degree} (or \emph{weight}) $\omega(\mu)$ of $\mu\in\mathcal{C}$ and we set $\mathcal C_n=\{\mu\in\mathcal C\mid\omega(\mu)=n\}$. If $\card(\omega^{-1}(0))=1$ then we denote by $\epsilon$ the unique element of $\mathcal{C}$ of weight $0$. We also set $\mathcal{C}^+:=\mathcal{C}\setminus \{\epsilon\}$.\\ \\
 
\subsection{Labelled combinatorial classes}
Recall that the symmetric group $\mathcal S_n$ is the group of bijections of $\{1,\dots,n\}$. It is a group of order $n!$ whose each element is denoted by the word of its images. For instance, the cycle sending $1$ to $2$, $2$ to $3$, and $3$ to $1$ is denoted by $231$. Obviously, the set of permutations is closed by composition and all its elements are invertible.\\
 Formally, a \emph{labelled combinatorial class} is a combinatorial class endowed with a sequence  $(\rho_n)_{n\in\mathbb N}$ of representations $\rho_n$ of the symmetric group $\mathcal S_n$ (\emph{i.e.} an application associating a map $\rho_n(\sigma):\mathcal C_n\rightarrow\mathcal C_n$ to each permutation $\sigma\in\mathcal S_n$ in such a way that $\rho_n(\sigma\circ\sigma')=\rho_n(\sigma)\circ\rho_n(\sigma')$).  An equivalent way (see \eg \cite{Flajolet}) to define labelled combinatorial class consists in considering that each element of $\mathcal C_n$ is a graph whose vertices are labelled by numbers from $1$ to $n$; the image of a permutation by the underlying representation is just the permutation of the labels.\\
 Let $(\mathcal C,\omega)$ and  $(\mathcal C',\omega')$ be two labeled combinatorial class. If the sets $\mathcal C$ and $\mathcal C'$ are disjoint, then we define the class $\mathcal C\boxplus \mathcal C'=(\mathcal C\cup\mathcal C',\omega'')$ with $\omega''(e)=\omega(e)$ if $e\in \mathcal C$ and $\omega(e)=\omega'(e)$ if $e\in\mathcal C'$. One extends to the case where $\mathcal C\cap\mathcal C'\neq \emptyset$ by replacing $\mathcal C'$ by a copy which is disjoint of $\mathcal C$ in the definition of $\mathcal C\boxplus \mathcal C'$.\\
 We also define $\mathcal C\boxtimes\mathcal C'$, i.e., the combinatorial class such that the elements of $(\mathcal C\boxtimes\mathcal C')_n$ are the pairs $(e,e')$ where $e$ is obtained by relabeling an element of $\mathcal C_i$ and $e'$ is obtained by relabeling an element of $\mathcal C'_j$, with $i+j=n$ such that the set of the labels in $(e,e')$ is $\{1,\dots,n\}$ and each relabeling preserves the initial order on the vertices. The degree of $(e,e')$ in $\mathcal C\boxtimes\mathcal C'$ is the sum of the degree of the respective preimage of $e$ and $e'$ in $\mathcal C$ and $\mathcal C'$.
 As a special case, for each labeled class $\mathcal C$, we denote  $\mathcal C^\bullet=\bullet\boxtimes\mathcal C$, where $\bullet$ is the class of the unique element of which, denoted also by $\bullet$, has degree $1$.
\subsection{The exponential generating function of a combinatorial class}

The exponential generating function (EGF) of a combinatorial class $\mathcal{C}$ is the exponential generating function of the numbers $C_n=\card(\mathcal{C}_n),$ in other words \[S_{\mathcal C}(x)=\sum\limits_{n\geq 0}C_n \dfrac{x^n}{n!}=\sum\limits_{\mu\in\mathcal{C}}\dfrac{x^{|\mu|}}{|\mu|!}.\]
We say that two classes are \emph{isomorphic} if their EGF are equal
\[
\mathcal C\equiv\mathcal C'\Leftrightarrow S_{\mathcal C}=S_{\mathcal C'}\Leftrightarrow \forall n\in\mathbb N, C_n=C'_n.
\]
Classically, we have \cite{Flajolet}
\[
S_{\mathcal C\boxplus\mathcal C'}=S_{\mathcal C}+S_{\mathcal C'}\mbox{ and }
S_{\mathcal C\boxtimes\mathcal C'}=S_{\mathcal C}S_{\mathcal C'}.
\]
\subsection{Labelled sequences}

%The basic building blocks from which many complex objects can be constructed are the three combinatorial classes $\Set,\Seq$ and $\Cyc.$ Recall that the exponential generating function for these classes are given by $\dfrac{1}{1-B(x)}$, $\exp{B(x)}$, and respectively, $\log\dfrac{1}{1-B(x)}$ (see \cite{Flajolet}).

If $\mathcal C^+$ is a labelled combinatorial class such that $\mathcal C^+_0=\emptyset$ then we define, up to an equivalence, the class $\Seq(\mathcal C^+)$ of labelled sequences by the equation
\[
\Seq(\mathcal C^+)\equiv [\ ]\boxplus (\mathcal C^+\boxtimes \Seq(\mathcal C^+)),
\]
where  $[\ ]$ denotes the class having a single element $\epsilon$ which is degree $0$.\\
It is easy to show that such a class exists and that its associated exponential generating function is
\[
S_{\Seq(\mathcal C^+)}(x)=\frac1{1-S_{\mathcal C^+}(x)}.
\]
From a combinatorial point of view, the elements of $\Seq(\mathcal C^+)_n$ are $k$-tuple $[\mu_1,\dots,\mu_k]$ where each $\mu_i$ is obtained by an order preserving relabelling of an element of $\mathcal C_{j_i}^+$ ,in such a way that the whole set of labels in $[\mu_1,\dots,\mu_k]$ is $\{1,\dots,n\}$ (as a consequence one has $\sum_{i=1}^k j_i=n$). So to any element $s=[\mu_1,\dots,\mu_k]\in \Seq(\mathcal C^+)_n$ we associate an ordered partition $\Pi=[\Pi_1,\dots,\Pi_k]$ of size $n$ such that each $\Pi_i$ is the set of the labels of $\mu_i$.

Let us be a little more precise. A \emph{labelled} list of elements of $\mathcal{C}^+$ is a list $L=[\mu_1,\mu_2,\ldots,\mu_k]\in\Seq(\mathcal{C}^+)$ with each $\mu_i$ associated to a set $\Omega_i$ of non-negative integers with ${\rm card}(\Omega_i)=|\mu_i|$ and $\Omega_i\cap\Omega_j = \emptyset$ for all $i,j$. A \emph{standard} labelled list of elements of $\mathcal{C}^+$ is a labelled list $L$ of elements of $\mathcal{C}^+$ such that the set of all the labels of $L$ is $\{1,\ldots,|L|\}$.

\subsection{Labelled sets and labelled cycles}
We consider the labelled combinatorial classe $\Set(\mathcal C^+)$ such that $\Set(\mathcal C^+)_n$ is the quotient of the set $\Seq(\mathcal C^+)_n$ by the relation $[\mu_1,\dots,\mu_k]\equiv_S [\mu_{\sigma(1)},\dots,\mu_{\sigma(k)}]$ for any permutation $\sigma\in\mathcal S_k$. Straightforwardly, each element of $\Set(\mathcal C)$ can be represented by a set of (order preserved) relabelled elements of $\mathcal C^+$. The exponential generating function
\[
S_{\Set(\mathcal C^+)}(x)=\exp\{S_{\mathcal C^+}(x)\}.
\] is  easily deduced from the construction.\\
If one consider the equivalence relation generated by $[\mu_1,\dots,\mu_k]\equiv_C [\mu_{2},\dots,\mu_{k},\mu_1]$, the one obtains an other labelled combinatorial class $\Cyc(\mathcal C^+)$ whose elements can be represented by necklace of (order preserved) relabelled elements of $\mathcal C^+$. We denote a necklace by $(\mu_1,\dots,\mu_k)=(\mu_2,\dots,\mu_k,\mu_1)$. Again, the generating series 
\[
S_{\Cyc(\mathcal C^+)}(x)=\log\left\{\frac1{1-S_{\mathcal C^+}(x)}\right\}.
\]
is deduced from the construction.  
\section{Set partitions and related constructions}
%{\tt ceci est utile pour comprendre la derniere section}
\subsection{Three constructions based on set partitions}
A \emph{set partition} of size $n$ is a set $\pi=\{\pi_1,\dots,\pi_k\}$ such that $\pi_1\cup\cdots\cup\pi_k=\{1,\dots,n\}$ and $\pi_i\cap\pi_j=\emptyset$ for any two indices $1\leq i\neq j\leq k$.
The set of set partitions $\mathcal S_p$ endowed with the size is a combinatorial classes satisfying
$   \mathcal S_p\equiv\Set(\Set(\bullet)^+)
$, and so, $S_{\mathcal S_p}(x)=\exp(\exp(x)-1)$. The numbers $\mathtt Sp_n= 	1, 1, 2, 5, 15, 52, 203, 877, 4140\dots$ are the well known Bell numbers (see sequence A000110 in \cite{Sloane}).\\
If $\mathcal C^+$ is a labeled combinatorial class such that $\mathcal C^+_0=0$ and $\mathcal P\equiv \Set(\mathcal C^+)$ then the definitions above allows to associate to each element $p=\{p_1,\dots,p_k\}\in \mathcal P_n$ a set partitions $\pi(p)=\{\labels(p_1),\dots,\labels(p_k)\}$ of size $n$  where for each $1\leq i\leq k$, $\labels(p_i)$ denotes the set of the labels of $p_i$.

An \emph{ordered partition} of size $n$ is a sequence $\Pi=[\Pi_1,\dots,\Pi_k]$ of non empty sets such that $\{\Pi_1,\dots,\Pi_k\}$ is a set partition of $\{1,\dots,n\}$. The set of ordered partitions $\mathcal{O}p$ endowed with the size is a combinatorial classes satisfying
$
\mathcal{O}p\equiv\Seq(\Set(\bullet)^+)$ and $S_{\mathcal{O}p}(x)=\frac1{2-\exp(x)}$.

The numbers $\mathtt Op_n=1, 1, 3, 13, 75, 541, 4683, 47293, 545835,\dots$ are the Fubini numbers (see sequence 
A000670 in \cite{Sloane}). If  $\mathcal L\equiv \Seq(\mathcal C^+)$ then the definitions above allows to associate to each element $\ell=[\ell_1,\dots,\ell_k]\in \mathcal L_n$ a set partition $\Pi(\ell)=[\labels(\ell_1),\dots,\labels(\ell_k)]$ of size $n$.

A \emph{cyclic partition} of size $n$ is a necklace $\mathfrak p=(\mathfrak p_1,\dots,\mathfrak p_k)$ such that $\{\mathfrak p_1,\dots,\mathfrak p_k\}$ is a set partition of size $n$. The set of ordered partitions $\mathcal{C}p$ endowed with the size is a combinatorial classes satisfying
$
\mathcal{C}p\equiv\Cyc(\Set(\bullet)^+)$ and $S_{\mathcal{C}p}(x)=\log(\frac1{2-\exp(x)})$.

The numbers $\mathtt Cp_n=1, 1, 3, 13, 75, 541, 4683, 47293, 545835,\dots$ are listed in sequence 
A000670 \cite{Sloane}. If  $\mathcal Ne\equiv \Cyc(\mathcal C^+)$ then the above definitions allow us  to associate to each element $c=(c_1,\dots,c_k)\in \mathcal (Ne)_n$ a cyclic partition $\mathfrak p(c)=(\labels(c_1),\dots,\labels(c_k))$ of size $n$.
\subsection{An explicit isomorphism}
%{\tt J'ai (re)changé les max en min dans les définitions de cette sous-section pour que cela colle avec les mots de Lyndon}

From the generating series we have $\Set(\mathcal Cp)\equiv\mathcal Op$. Indeed, this equality translates in terms of generating function as $\exp\left(\log\left(\frac1{2-e^x}\right)\right)=\frac1{2-e^x}$. In order to understand a more general identity introduced later in the paper, we make explicit this bijection. 
Assume that $c=\{c^{(1)},\dots,c^{(k)}\}\in \Set(\mathcal Cp)$. If $c^{(i)}=(c^{(i)}_1,\dots,c^{(i)}_{h_i})$ then we consider  $\sigma_i$  the only circular permutation on the indices $\{1,\dots,h_i\}$ such that $\min\bigcup_j\mathtt{labels}(c^{(i)}_{(j)})=\min\mathtt{labels}(c^{(i)}_{\sigma_i^{-1}(1)})$. In other words, if $\ell_i=[\ell^{(i)}_1,\dots, \ell^{(i)}_{h_i}]=[c^{(i)}_{\sigma_i(1)},\dots,c^{(i)}_{\sigma_i(h_i)}]$ then $\min\bigcup_j\mathtt{labels}(\ell^{(i)}_{(j)})=\min\mathtt{labels}(\ell^{(i)}_{1})$. Now, consider the unique permutation $\rho\in\mathcal S_k$ such that \[\min\mathtt{labels}(c^{(\rho^{-1}(1))})> \min\mathtt{labels}(c^{(\rho^{-1}(2))})>\cdots> \min\mathtt{labels}(c^{(\rho^{-1}(k))})\] and set
\begin{equation}\mathtt{stol}(c)=[\ell^{(\rho(1))}_1,\dots,\ell^{(\rho(1))}_{h_{\rho(1)}},\dots,\ell^{(\rho(k))}_k,\dots,\ell^{(\rho(k))}_{h_{\rho(k)}}]\in\mathcal Op.\end{equation}
For instance,  
\begin{eqnarray*}\mathtt{stol}(\{(\{11\},\{2,5\},\{10\}),(\{6\},\{1,3,4\},\{7,9\}),(\{8,12\})\})=\\\ [\{8,12\},\{2,5\},\{10\},\{11\},\{1,3,4\},\{7,9\},\{6\}].\end{eqnarray*}
Let $\ell=[\ell_1,\dots,\ell_k]\in\mathcal Op$ and $1=i_0\leq \dots\leq i_{h-1}<i_{h}=k+1\in\{1,\dots,k+1\}$ be the set of indices satisfying 
\begin{equation}\label{mineq}\min\bigcup_{i< i_{j+1}}\mathtt{\labels}(\ell_i)>
\min\bigcup_{i\geq i_j}\mathtt{\labels}(\ell_i)\end{equation} with $h$ maximal.\\ For instance, the indices associated to $[\{8,12\},\{2,5\},\{10\},\{11\},\{1,3,4\},\{7,9\},\{6\}]$ are $1\leq 2\leq 5< 8$. We define 
\begin{equation}
    \mathtt{ltos}(\ell)=\{c_1,\dots,c_k\}\in\Set(\mathcal Cp),\end{equation}
where $c_j$ denotes the necklace $(\ell_{i_{j-1}},\dots,\ell_{i_{j}-1})$.
For instance,  
\begin{eqnarray*}\mathtt{ltos}([\{8,12\},\{2,5\},\{10\},\{11\},\{1,3,4\},\{7,9\},\{6\}])=\\\ \{(\{1,3,4\},\{7,9\},\{6\}),(\{2,5\},\{10\},\{11\}),(\{8,12\})\}.\end{eqnarray*}
It is easy to check that $\mathtt{ltos}(\mathtt{stol}(c))=c$ and  $\mathtt{stol}(\mathtt{ltos}(\ell))=\ell$. So we have
\begin{proposition}
The map $\mathtt{stol}$ is an isomorphism of combinatorial classes and  $\mathtt{ltos}$ is its reverse map.
\end{proposition}
\subsection{About Lyndon words}
Recall that the free monoid (\eg\cite{Loth}) $\Sigma^*$ on a set $\Sigma$ is the monoid whose elements are all the finite sequences endowed with the catenation product $\cdot$ that consists of pasting  one sequence to the right of another. The empty sequence plays the role of the identity element. For instance, in the free monoid $\{a,b\}^*$ we have $[a,b,a,a,b]\cdot[b,a,a]=[a,b,a,a,b,b,a,a]$. In literature,  brackets and commas are often omitted; the elements of a free monoid are then noted as juxtapositions of letters called words (the empty word, noted by $\varepsilon$, corresponds to the sequence $[\ ]$). The name of the free monoid comes from the fact that it fulfills the universal property, that is every monoid having a generating set in bijection with $\Sigma$ is isomorphic to a quotient of $\Sigma^*$.\\
Any pair of sequences under the form $u\cdot v$ and $v\cdot u$ are said conjugate. In other words, the conjugates of a sequence are all its circular shift. This is obviously an equivalence relation that preserves the periods, \ie the conjugate sequences of $u^{\cdot k}$ are exactly the sequences $v^{\cdot k}$ where $v$ is conjugate to $u$. In terms of combinatorial class the free monoid is nothing but $\Seq(\Sigma)$ and its quotient by conjugation is $\Cyc(\Sigma)$.\\
Assume that the alphabet $\Sigma$ is totally ordered by the order $<$. Then the free monoid  is totaly ordered with the lexicographic order $\prec$. The minimal element for the lexicographic order is the empty sequence $[\ ]$ and we have $[a]\cdot u\prec[b]\cdot v$ if $a<b$ or $a=b$ and $u\prec v$.\\
A \emph{Lyndon words} a non periodic sequence which is minimal in its conjugacy class. Their name comes from the mathematician Roger Lyndon who studied them in 1954\cite{Lyndon1}. Nevertheless, it should be noted that they had been introduced a year earlier by  Anatoly Shirshov \cite{Shir1}. Lyndon words play a very important role for understanding of free groups \cite{Lyndon2}, free associative algebras, and free Lie algebras \cite{Shir2}. Readers may refer to \cite{Reut} for a rather complete survey.\\
Among all the properties of Lyndon's words, one of the most interesting is that they play for the free monoid the same role as prime numbers play for integers. This property is that any sequence factorizes as a unique weakly decreasing catenation of  Lyndon words \cite{Schutz}. In other words, the free monoid $\Sigma^*$ is in bijection with the multisets of aperiodic sequences over $\sigma$. For instance, if we assume $a<b$ the sequence $
u=[a,b,a,b,b,a,b,a,b,a,a,a,b,a,b,a]$ factorizes as $u=[a,b,a,b,b]\cdot[a,b]\cdot[a,b]\cdot[a,a,a,b,a,b]\cdot [a]$. This means that the sequence $u$ is assimilated to the multiset \\ $\{(a,a,a,b,a,b),(a,b,a,b,b),(a,b),(a,b),(a)\}$ (remark the multiplicity of $(a,b)$).
It is interesting to note that this correspondence is precisely the one that is calculated when applying $\mathtt{ltos}$. Indeed, let $\ell=[\ell_1,\dots,\ell_k]\in\mathcal Op$, the alphabet  $\Sigma=\{\ell_1,\dots,\ell_k\}$ is totally ordered by $\ell_i<\ell_j$ if and only if $\min\ell_i<\min\ell_j$.  In fact, since each numbers of $\{1,\dots,n\}$ appears only one time in the sequence, only the minimal elements the sets are relevant and all works as if our alphabet be $\{1,\dots,n\}$. For instance, $[\{8,12\},\{2,5\},\{10\},\{11\},\{1,3,4\},\{7,9\},\{6\}]$ is assimilated to $[8,2,10,11,1,7,6]$.
The indices of equation (\ref{mineq}), except the larger which is not relevant,  indicate where to catenate in order to apply the complete factorization. In our example, we found the indices $\{1,2,5,8\}$ and, then, we have  $[8,2,10,11,1,7,6]=[8]\cdot[2,10,11]\cdot[1,7,6]$. Notice that, since the components are two by two distinct, the bijection with multisets of cycles class sends the sequences we consider on set of necklaces. For instance, $[8]\cdot[2,10,11]\cdot[1,7,6]\sim \{(1,7,6),(2,10,11),(8)\}$. We recover the $\mathtt{ltos}(\ell)$ by replacing each integer by the set of which it is the minimum. In our example we have $$\{(1,7,6),(2,10,11),(8)\}\rightarrow 
\{(\{1,3,4\},\{7,9\},\{6\}),(\{2,5\},\{10\},\{11\}),(\{8,12\})\}.$$
Of course, this may seem like a very sophisticated way to revisit the bijection of the previous section. Nevertheless, this remark is valuable because it will allow us to link our constructions to notions of algebras (enveloping algebras, Hopf algebras, Lie algebras of primitive elements etc.) that we will explore in future works.
\subsection{Set of cycles and sequences}\label{sec:setcyc-seq}
Let $\mathcal C^+$ be a labeled combiatorial sequences such that $\mathcal C^+_0=\emptyset$. We define
$    \mathcal J=\Set(\Cyc(\mathcal C^+))$ and $\mathcal S=\Seq(\mathcal C^+)$.

Let us show that the map $\mathtt{stol}$ allows us to compute an explicit isomorphism from $\mathcal J$ to $\mathcal S$. We define $\mathtt{jtoset}_{\mathcal C^+}:\mathcal J\rightarrow \Set(\mathcal Cp) $ by \begin{eqnarray*}
\mathtt{jtoset}_{\mathcal C^+}(\{(c^{(1)}_1,\dots,c^{(1)}_{h_1}),\cdots,(c^{(k)}_1,\dots,c^{(k)}_{h_k})\})=\\
\{(\mathtt{labels}(c^{(1)}_1),\dots,\mathtt{labels}(c^{(1)}_{h_1})),\cdots,(\mathtt{labels}(c^{(k)}_1),\dots,\mathtt{labels}(c^{(k)}_{h_k}))\}.\nonumber
\end{eqnarray*}
We define also $\mathtt{jtoseq}_{\mathcal C^+}:\mathcal J\rightarrow\mathcal S$ such that $$\mathtt{jtoseq}_{\mathcal C^+}(\{(c^{(1)}_1,\dots,c^{(1)}_{h_1}),\cdots,(c^{(k)}_1,\dots,c^{(k)}_{h_k})\})=\ell$$ is the unique permutation of the vector $[c^{(1)}_1,\dots,c^{(1)}_{h_1},\cdots,c^{(k)}_1,\dots,c^{(k)}_{h_k}]$ such that $\Pi(\ell)=\mathtt{stol}(\mathtt{jtoset}_{\mathcal C^+}((\{(c^{(1)}_1,\dots,c^{(1)}_{h_1}),\cdots,(c^{(k)}_1,\dots,c^{(k)}_{h_k})\})))$.

 Let $c=\{(c_1^{(1)},c_2^{(1)}),(c_1^{(2)}),(c_1^{(3)},c_2^{(3)})\}$ be a set of labelled graphs as shown in Fig. \ref{fig:labeled-graph}    
\begin{figure}[!ht]
 	\begin{center}
 	\resizebox{\textwidth}{!}{
    \begin{tikzpicture}[thin,sibling distance=10mm,
      every circle node/.style={minimum size=1.5mm,inner sep=0mm}, level distance=10mm]
 
\node[circle,fill,label=left:$3$] at (0,0){};
\node[circle,fill,label=below:$7$] at (1,0){};
\node[circle,fill,label=right:$5$] at (2,0){};
\node[circle,fill,label=left:$20$] at (1,1){};
\draw (0,0) -- (1,0) -- (2,0); 
\draw (1,0) -- (1,1);

\node[circle,fill,label=left:$9$] at (4,0){};
\node[circle,fill,label=right:$14$] at (6,0){};
\node[circle,fill,label=left:$17$] at (5,1){};
\node[circle,fill,label=left:$11$] at (4,2){};
\node[circle,fill,label=right:$19$] at (6,2){};
\draw (4,0) -- (6,0) -- (5,1)-- (4,2);
\draw (4,0) -- (5,1) -- (6,2);

\node[circle,fill,label=left:$1$] at (8,0){};
\node[circle,fill,label=left:$18$] at (8,1){};
\draw (8,0) -- (8,1); 

\node[circle,fill,label=left:$2$] at (10,0){};
\node[circle,fill,label=right:$8$] at (12,0){};
\node[circle,fill,label=left:$13$] at (11,1){};
\draw (10,0) -- (12,0) -- (11,1)-- (10,0);

\node[circle,fill,label=left:$4$] at (14,0){};
\node[circle,fill,label=right:$16$] at (16,0){};
\node[circle,fill,label=left:$6$] at (14,1){};
\node[circle,fill,label=right:$12$] at (16,1){};
\node[circle,fill,label=below:$15$] at (15,1) {};
\node[circle,fill,label=left:$10$] at (15,2){};

\draw (15,2) -- (15,1) -- (14,1)-- (14,0)-- (16,0) -- (16,1)-- (15,1);

\node[thick] at (1,-1) {$\Large{c_1^{(1)}}$};
\node[thick] at (5,-1) {$\Large{c_2^{(1)}}$};
\node[thick] at (8,-1) {$\Large{c_1^{(2)}}$};
\node[thick] at (11,-1) {$\Large{c_1^{(3)}}$};
\node[thick] at (15,-1) {$\Large{c_2^{(3)}}$};
          \end{tikzpicture}
          }
         % }
         \caption{Five labelled graphs}\label{fig:labeled-graph}
    \end{center}
\end{figure}
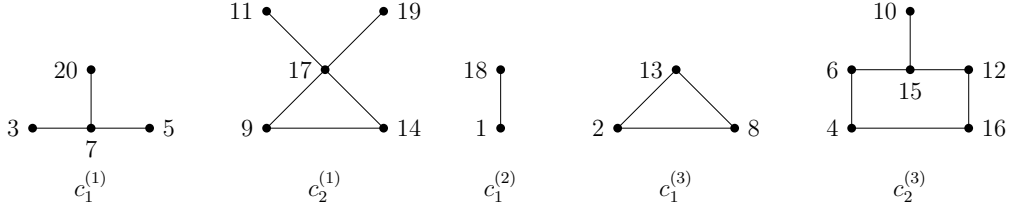

Then $\mathtt{jtoset}_{\mathcal{C^+}}(c)$ equals \begin{eqnarray*}
\{(\{3,7,5,20\},\{9,14,17,19,11\}),(\{1,18\}),(\{2,8,13\},\{4,16,12,15,10,6\})\}.
\end{eqnarray*}

Furthermore, $\Pi(l)=\mathtt{stol}\left(\mathtt{jtoset}_{\mathcal{C^+}}(c)\right)=$
\begin{eqnarray*}
[\{3,7,5,20\},\{9,14,17,19,11\},\{2,8,13\},\{4,16,12,15,10,6\},\{1,18\}].
\end{eqnarray*}
Since $\mathtt{stol}$ is one to one, the equality on generating functions allows us to deduce that $\mathtt{jtoseq}_{\mathcal C^+}$ is an isomorphism of combinatorial classes. The inverse bijection $\mathtt{seqtoj}_{\mathcal C^+}:\Set(\mathcal Cp)\rightarrow \mathcal J$ is defined by $$\mathtt{seqtoj}_{\mathcal C^+}(\ell)=\{(c^{(1)}_1,\dots,c^{(1)}_{h_1}),\cdots,(c^{(k)}_1,\dots,c^{(k)}_{h_k})\}$$ where $[c^{(1)}_1,\dots,c^{(1)}_{h_1},\cdots,c^{(k)}_1,\dots,c^{(k)}_{h_k}]$ is the unique permutation of $\ell$ such that
$\mathtt{ltos}(\Pi(\ell))=\mathtt{jtoset}_{\mathcal C^+}(\{(c^{(1)}_1,\dots,c^{(1)}_{h_1}),\cdots,(c^{(k)}_1,\dots,c^{(k)}_{h_k})\})$.

In the aforementioned example, we deduce the indices of the minimum elements in $\Pi(l)$ being $1<3<5$. Hence,
\begin{align*}
 \mathtt{ltos}(\Pi(l)) &=\{(\{1,18\}),(\{2,8,13\},\{4,16,12,15,10,6\}),(\{3,7,5,20\},\{9,14,17,19,11\})\}\\
   						&=\mathtt{jtoset}_{\mathcal{C^+}}(c).
\end{align*}

We summarize the results of this section in the following theorem
\begin{theorem}\label{th1}  The maps which make commuting the following diagram are explicit isomorphisms of combinatorial classes 
\begin{equation}
 \Set(\Cyc(\mathcal C^+))\mathop{\rightleftarrows}^{\mathtt{jtoseq}_\mathcal C^+}_{\mathtt{seqtoj}_{\mathcal C^+}}\Seq(\mathcal C^+).
\end{equation}
\end{theorem}
%\textbf{On devrait plutot avoir SEQ(C) au lieu de SET(C)!}

\section{Labelled and unlabelled trees}\label{sec:trees}

We illustrate the previous result by investigating the combinatorial classes $\mathcal{R}$ satisfying \begin{equation}\Set(\Cyc(\mathcal R^\bullet))\equiv \Seq(\mathcal R^\bullet)\equiv \mathcal{R}.\label{eq:Cat-eq}\end{equation}
This isomorphism can be translated into the following functional equation:\begin{equation}\label{functeqR} \exp\left\{\log\left\{\frac1{1-xS_{\mathcal R}(x)}\right\}\right\}=\dfrac{1}{1-xS_{\mathcal R}(x)}=S_{\mathcal R}(x).\end{equation} 
This equation has a unique solution 
 \begin{equation}S_{\mathcal R}(x)=\dfrac{1-\sqrt{1-4x}}{2x}\end{equation} 
which is also the ordinary generating function of the Catalan numbers $\mathtt C_n=\frac1{n+1}\binom{2n}n$  \cite{Catalan,Flajolet}. Hence, $\mathcal R$ is unique up to an isomorphism and
\begin{equation}
R_n=\dfrac{(2n)!}{(n+1)!}.    
\end{equation}
The sequence of $R_n$ is  
\[1, 1, 4, 30, 336, 5040, 95040, 2162160  \quad \href{http://oeis.org/A001761}{A001761} \cite{Sloane} .\]

\subsection{Labelled trees from unlabelled trees}
\begin{definition}
A \emph{tree} is a list of trees (possibly empty) connected to a node, called its root, by an edge (also called branch). Notice that this is a valid recursive definition which base case is a root together with an empty list. The \emph{degree} $\omega(t)$ of a tree $t$ is the number of its edges or, equivalently the number of its nodes which are not its root. 
\end{definition}
Let $\D$ be the set of trees. There are a finite number of trees having a given degree, so the pair $(\D,\omega)$ is a (unlabelled) combinatorial class.
The number $D_n$ is known to be the Catalan number $\mathtt C_n$ (see \eg \cite{Catalan}). So the ordinary generating function of the class $\mathcal D$, \ie
\begin{equation}
    S_{\mathcal D}^{ord}(x)=\sum_{n\geq 0}D_nx^n,
\end{equation}
fulfills the same functional equation (\ref{functeqR}) as the exponential generating function  of $\mathcal R$. So each $\mathcal R_n$ is in one to one correspondence with $\mathcal D_n\times\mathcal S_n$. This suggests that one can exhibit an explicit realization of the class $\mathcal R$ by labeling the nodes which are not the root of each tree $t\in\mathcal D$ by $\{1,\dots,\omega(t)\}$, without repetition and in any possible way.

\subsection{Shifted structure}
The class $\mathcal R^\bullet$ is isomorphic to the class $\mathcal R_r$ of trees with labeled root endowed with the weight $\omega_r$ counting the total number of nodes, including the root.
%As the root is a particular node it is natural to associate it with a label.
More precisely, for a given $n$, a tree of $(\R_r)_{n+1}$ is obtained by labelling the root of a tree in $\R_n$ with any of the possible value from the set $\{1,\dots,n+1\}$ and relabel, if necessary, the nodes with respect to the order induced by the initial permutation. 
\begin{example}
Let $t$ be a rooted labelled tree in $\R_5$, and its associated permutation is $\pi=(1,2,4,5,3)$. Then the set of all rooted labeled trees where the root is labeled obtained from $t$ is given by the set of permutations
$\{(6,1,2,4,5,3),\\(5,1,2,4,6,3),(4,1,2,5,6,3),(3,1,2,5,6,4),(2,1,3,5,6,4),(1,2,3,5,6,4)\}.$
\end{example}

In terms of generating function, this operation leads to 
\[S_{\R_r}(x)=S_{\R^\bullet}(x)=xS_{\R}(x)=\dfrac{1-\sqrt{1-4x}}{2}\] and so \[(R_r)_n=n!C_{n-1}=\dfrac{(2n-2)!}{(n-1)!},\] for any  $n\ge 1.$ The sequence of ${R_r}_n$ is given by 

\[1, 2, 12, 120, 1680, 30240  \quad \href{http://oeis.org/A001761}{A001813}\ \cite{Sloane} .\] 

We insists on the fact that at this point ${(\R_r)}_n$ counts trees having $n$ nodes including the root. 

\subsection{Hanging trees in necklaces}
A \emph{labelled necklace of planar trees} is
a necklaces on which trees are hung and all the nodes (comprising the roots) are labeled by $\{1,\dots,n\}$ where $n$ is the total numbers of nodes (comprising roots). We denote by $\mathcal N$ the set of such necklaces. The cyclic structure comes from the fact that a necklace is invariant  by rotation. The weight $\omega_N(\mathtt n)$ of a necklace $\mathtt n$ is the total number of the nodes, comprising roots, of the trees it contains. The pair $(\mathcal N,\omega_N)$ is a labeled combinatorial class that satisfies
\begin{equation}\label{N2R}
    \mathcal N\equiv\Cyc(\mathcal R^\bullet).
\end{equation}

We depict in Fig. \ref{fig:planar} elements of $\mathcal N$ that are equivalent under cyclic rotation. 

\begin{figure}[!ht]
 	\begin{center}
    \small
    \begin{tikzpicture}[thin,sibling distance=10mm,
      every circle node/.style={minimum size=1.5mm,inner sep=0mm}, level distance=10mm]
      
      \node[draw=none] at (-3,6) (root1) {}
        child { node [circle,fill,label=below:$1$] {}
          edge from parent
          };
    \node[draw=none] at (-2,5.1) (root1) {}      
        child { node [circle,fill,label=below:$5$] {}
          edge from parent
          };
    \node[draw=none] at (-1,6) (root1) {}      
          child { node [circle,fill,label=right:$2$] {}
           edge from parent
            child { node [circle,fill,label=below:$4$] {}
           edge from parent
             }
              child { node [circle,fill,label=below:$3$] {}
           edge from parent
             }
             };
             %3rd tree
           \node[draw=none] at (5,6) (root2) {}
        child { node [circle,fill,label=left:$2$] {}
          edge from parent
           child { node [circle,fill,label=below:$4$] {}
           edge from parent
             }
              child { node [circle,fill,label=below:$3$] {}
           edge from parent
             }
          };
           \node[draw=none] at (6,5.1) (root2) {}
        child { node [circle,fill,label=below:$1$] {}
           edge from parent
             };
           \node[draw=none] at (7,6) (root2) {}
        child { node [circle,fill,label=below:$5$] {}
           edge from parent
             };     
   %2nd tree     
           \node[draw=none] at (1,6) (root3) {}
        child { node [circle,fill,label=below:$5$] {}
           edge from parent
             };
             \node[draw=none] at (2,5.1) (root3) {}
        child { node [circle,fill,label=left:$2$] {}
          edge from parent
           child { node [circle,fill,label=below:$4$] {}
           edge from parent
             }
              child { node [circle,fill,label=below:$3$] {}
           edge from parent
             }
          };
          \node[draw=none] at (3,6) (root3) {}
        child { node [circle,fill,label=below:$1$] {}
           edge from parent
             };     
                            
         \draw (-2,6) circle (1cm);    
\draw (2,6) circle (1cm);    
\draw (6,6) circle (1cm);    
          \end{tikzpicture}
          \caption{Three rooted labeled trees equivalent under cyclic rotation}\label{fig:planar}
    \end{center}
\end{figure}
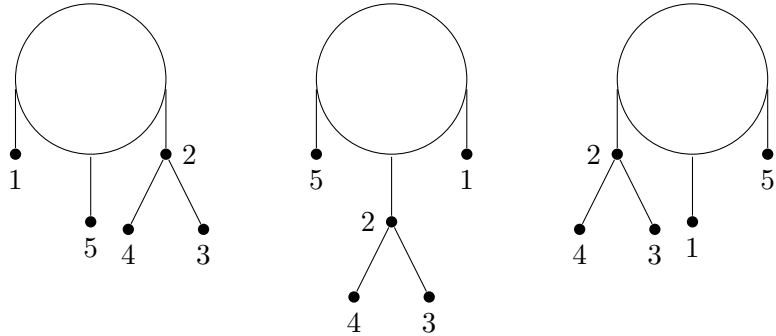

Notices that Labeled necklaces of rooted trees appear under the name of "planar labelled trees" in the work of Miloudi \cite{Miloudi}. To be more precise, he studied combinatorial class which is straightforwardly isomorphic to $\mathcal N$ and he proved  $N_n=(2n-3)!/(n-1)!$ for any $n\geq 2$ and $N_1=1.$  We recover this result from the interpretation of (\ref{N2R}) in terms of generating function.  Indeed, \[S_{\Cyc(\mathcal R^{\bullet})}(x)=\log\dfrac{1}{1-\dfrac{1-\sqrt{1-4x}}{2}}=\log\dfrac{2}{1+\sqrt{1-4x}}=\log{\dfrac{1-\sqrt{1-4x}}{2x}.}\]Hence, 
\begin{equation}\label{eq:cyc-log}
    S_{\Cyc(\mathcal R^{\bullet})}(x)=\log \mathcal R(x),
\end{equation}
and the exact formula for $N_n$ is obtained by expanding the function as a Taylor series.\\
These numbers are also mentioned by Wolfdieter Lang in \cite{Sloane}, see the sequence below
 \[1, 1, 3, 20, 210, 3024, 55440, 1235520, 32432400\quad \href{http://oeis.org/A006963}{A006963}.\]

%For example there are $20$ rooted planar labelled trees with $3$ labelled nodes and one root, as given by $P_4$. 
\medskip
%and it is given by the sequence \[1, 1, 3, 20, 210, 3024, 55440, 1235520, 32432400\quad \href{http://oeis.org/A006963}{A006963}.\]

%\subsection{Counting Labelled Trees}

%The number of labelled trees with $n$ nodes, also known as Cayle's formula, equals $L_n=n^{n-2}$ and it is given by the sequence \[1, 1, 3, 16, 125, 1296, 16807, 262144\quad \href{http://oeis.org/A000272}{A000272}.\]

%Since at each node we can add an extra edge joining the node with the root we have that the number of rooted labelled trees with $n$ nodes equals $R_n^*=n^{n-1}$ and it is given by the sequence  \[1,2,9,64,625,7776, 117649, 2097152\quad \href{http://oeis.org/A000169}{A000169}.\] 

We now have all the material to make explicit the isomorphisms suggested by (\ref{eq:Cat-eq}). To this aim we consider jewellery boxes which are sets of necklaces and forests which are sequences of trees. More formally, in terms of combinatorial classes we define
$\mathcal J=\Set(\mathcal N)$ and $\mathcal F=\Seq(\mathcal R^\bullet)$.
Let $\mathtt{jtof}=\mathtt{jtoseq}_{\mathcal R^\bullet}:\mathcal J\rightarrow\mathcal F$ and its inverse bijection $\mathtt{ftoj}=\mathtt{seqtoj}_{\mathcal R^\bullet}:\mathcal F\rightarrow \mathcal J$.
An explicit isomorphism $\mathtt{rtof}:\mathcal R\rightarrow \mathcal{F}$ is obtained by removing the root to any tree in $\mathcal R$. The reciprocal isomorphism $\mathtt{ftor}:\mathcal F\rightarrow \mathcal{R}$ consists in connecting all the trees of a given sequence to an additional node called the root.\\
All these constructions are summarized in the following result which is a corollary of Theorem \ref{th1}.
%
%Notice that the labelled trees (with respect to $\pi$) in the forest are ordered, order induced by $\min$ on $\pi.$ 
%
%
\begin{corollary}
  The maps which make commuting the following diagram are explicit isomorphisms of combinatorial classes 
\begin{equation}
 \mathcal J\mathop{\rightleftarrows}^{\mathtt{jtof}}_{\mathtt{ftoj}}\mathcal F\mathop{\rightleftarrows}^{\mathtt{ftor}}_{\mathtt{rtof}}\mathcal R.\label{eq:cor1}
\end{equation}
\end{corollary}

In Fig. \ref{fig:ex-jtof-ftoj} we illustrate the bijections from \eqref{eq:cor1} using two examples.

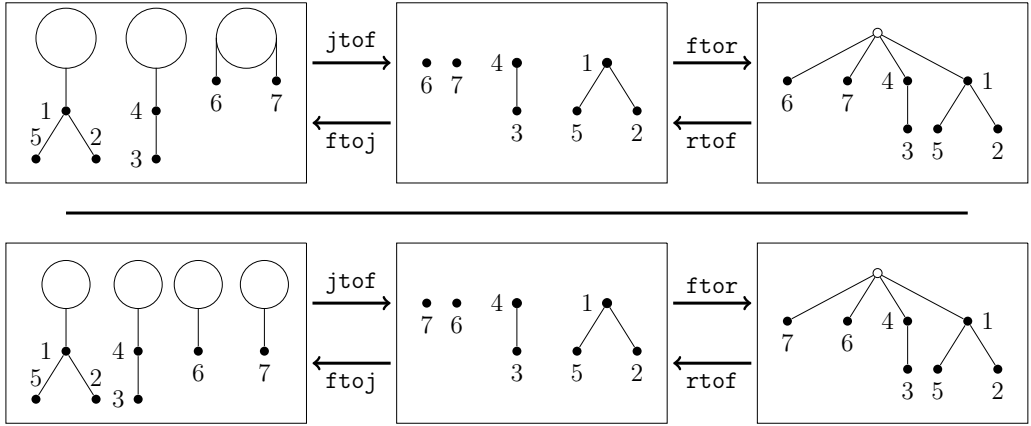
\begin{figure}[!ht]
 	\begin{center}
 	\resizebox{\textwidth}{!}{
    \begin{tikzpicture}[thin,sibling distance=10mm,
      every circle node/.style={minimum size=1.5mm,inner sep=0mm}, level distance=8mm]
           \draw (0,2.41) circle (0.5cm); 
          \node[circle] at (0.5,2.5) (root1) {}
        child { node [circle,fill,label=below:$7$] {}
           edge from parent
             };
             \node[circle] at (-0.5,2.5) (root1) {}
        child { node [circle,fill,label=below:$6$] {}
           edge from parent
             };
                   \draw (-1.5,2.41) circle (0.5cm); 
           \node[circle] at (-1-0.5,2) (root2) {}
        child { node [circle,fill,label=left:$4$] {}
          edge from parent
          child { node [circle,fill,label=left:$3$] {}
           edge from parent
             }
             };     
           \draw (-3,2.41) circle (0.5cm); 
           \node[circle] at (-3,2) (root3) {}
        child { node [circle,fill,label=left:$1$] {}
          edge from parent
        child { node [circle,fill,label=above:$5$] {}
           edge from parent
           }
        child { node [circle,fill,label=above:$2$] {}
         edge from parent
             }
            };     
           
           \draw (1,3) -- (-4,3) -- (-4,0) -- (1,0) -- (1,3); 
%second box            
          \node[circle,fill,label=below:$6$] at (-3+6,2) (root1) {};
        
          \node[circle,fill,label=below:$7$] at (-3+6.5,2) (root1) {};
        
           \node[circle,fill,draw,label=left:$4$] at (-1+6-0.5,2) (root2) {}
          child { node [circle,fill,label=below:$3$] {}
           edge from parent
             };     
           \node[circle,draw,fill,label=left:$1$] at (1+6-1,2) (root3) {}
        child { node [circle,fill,label=below:$5$] {}
           edge from parent
           }
        child { node [circle,fill,label=below:$2$] {}
         edge from parent
            };     
              \draw (1+6,3) -- (-4+6.5,3) -- (-4+6.5,0) -- (1+6,0) -- (1+6,3); 
    % 3rd box
           \node[circle,draw] at (10.5,2.5) (root3) {}
           child { node [circle,fill,label=below:$6$] {}
          edge from parent}
           child { node [circle,fill,label=below:$7$] {}
          edge from parent}
            child { node [circle,fill,label=left:$4$] {}
          edge from parent
          child { node [circle,fill,label=below:$3$] {}
           edge from parent
             }
             }
        child { node [circle,fill,label=right:$1$] {}
          edge from parent
        child { node [circle,fill,label=below:$5$] {}
           edge from parent
           }
        child { node [circle,fill,label=below:$2$] {}
         edge from parent
             }
            };     
                \draw (1+6+6,3) -- (-4+6+6+0.5,3) -- (-4+6+6+0.5,0) -- (1+6+6,0) -- (1+6+6,3); 

     \draw[->, line width=0.5mm] (1.1,2) -- node [above] {$\mathtt{jtof}$}(2.4,2);
     \draw[->, line width=0.5mm] (2.4,1) -- node [below] {$\mathtt{ftoj}$}(1.1,1);
     \draw[->, line width=0.5mm] (1.1+6,2) -- node [above] {$\mathtt{ftor}$}(2.4+6,2);
     \draw[->, line width=0.5mm] (2.4+6,1) -- node [below] {$\mathtt{rtof}$}(1.1+6,1);
     % end first exemple
     \draw[-,line width=0.5mm] (-3,-0.5)--(12,-0.5);
     %begin second example
      \draw (0.3,2.31-4) circle (0.4cm); 
          \node[circle] at (0.3,2-4) (root1) {}
        child { node [circle,fill,label=below:$7$] {}
           edge from parent
             };
      \draw (-0.8,2.31-4) circle (0.4cm); 
          \node[circle] at (-0.8,2-4) (root1) {}
        child { node [circle,fill,label=below:$6$] {}
           edge from parent
             };
      
                   \draw (-1.8,2.31-4) circle (0.4cm); 
           \node[circle] at (-1.8,2-4) (root2) {}
        child { node [circle,fill,label=left:$4$] {}
          edge from parent
          child { node [circle,fill,label=left:$3$] {}
           edge from parent
             }
             };     
           \draw (-3,2.31-4) circle (0.4cm); 
           \node[circle] at (-3,2-4) (root3) {}
        child { node [circle,fill,label=left:$1$] {}
          edge from parent
        child { node [circle,fill,label=above:$5$] {}
           edge from parent
           }
        child { node [circle,fill,label=above:$2$] {}
         edge from parent
             }
            };     
           
           \draw (1,3-4) -- (-4,3-4) -- (-4,0-4) -- (1,0-4) -- (1,3-4); 
%second box            
          \node[circle,fill,label=below:$7$] at (-3+6,2-4) (root1) {};
        
          \node[circle,fill,label=below:$6$] at (-3+6.5,2-4) (root1) {};
        
           \node[circle,fill,draw,label=left:$4$] at (-1+6-0.5,2-4) (root2) {}
          child { node [circle,fill,label=below:$3$] {}
           edge from parent
             };     
           \node[circle,draw,fill,label=left:$1$] at (1+6-1,2-4) (root3) {}
        child { node [circle,fill,label=below:$5$] {}
           edge from parent
           }
        child { node [circle,fill,label=below:$2$] {}
         edge from parent
            };     
              \draw (1+6,3-4) -- (-4+6.5,3-4) -- (-4+6.5,0-4) -- (1+6,0-4) -- (1+6,3-4); 
    % 3rd box
           \node[circle,draw] at (10.5,2.5-4) (root3) {}
           child { node [circle,fill,label=below:$7$] {}
          edge from parent}
           child { node [circle,fill,label=below:$6$] {}
          edge from parent}
            child { node [circle,fill,label=left:$4$] {}
          edge from parent
          child { node [circle,fill,label=below:$3$] {}
           edge from parent
             }
             }
        child { node [circle,fill,label=right:$1$] {}
          edge from parent
        child { node [circle,fill,label=below:$5$] {}
           edge from parent
           }
        child { node [circle,fill,label=below:$2$] {}
         edge from parent
             }
            };     
                \draw (1+6+6,3-4) -- (-4+6+6+0.5,3-4) -- (-4+6+6+0.5,0-4) -- (1+6+6,0-4) -- (1+6+6,3-4); 

     \draw[->, line width=0.5mm] (1.1,2-4) -- node [above] {$\mathtt{jtof}$}(2.4,2-4);
     \draw[->, line width=0.5mm] (2.4,1-4) -- node [below] {$\mathtt{ftoj}$}(1.1,1-4);
     \draw[->, line width=0.5mm] (1.1+6,2-4) -- node [above] {$\mathtt{ftor}$}(2.4+6,2-4);
     \draw[->, line width=0.5mm] (2.4+6,1-4) -- node [below] {$\mathtt{rtof}$}(1.1+6,1-4);
         
         \end{tikzpicture}}
         \label{fig:ex-jtof-ftoj}
\caption{A set of three labelled necklaces of planar trees (leftmost box) in bijection with a forest of labelled trees (middle box) in bijection with a rooted labelled tree (rightmost box). In the upper part of the figure the cycles are $\{(\{1,5,2\}),(\{3,4\}),(\{6\},\{7\})\}$, whereas in the lower part the cycles are $\{(\{1,5,2\}),(\{3,4\}),(\{6\}),(\{7\})\}$}
    \end{center}
\end{figure}

\subsection{Other recursive  tree-like combinatorial classes}\label{sec:rec-tree-classes}
The tree-like structure constructed from sequences is the most rigid one. It involves rooted trees that are embedded in a plan in such a way that the branches are always pointing downwards and so, the order of the sequences of the subtrees is relevant. If we relax the constraint of the orientation of the branches then the trees become invariant by rotation. In that context, a tree is a non-oriented graph without cycle with a privileged vertex called a root. Each root can be seen as a labeled necklace on which the subtrees are hanged, forming a sort of windmill (see Figure \ref{windmills} for examples).  In other words, a tree is either an isolated root or a root with a cycle of trees. The combinatorial class satisfies the isomorphism
\begin{equation}
    \mathcal W\equiv \bullet\boxplus\left(\bullet\boxtimes\Cyc(\mathcal W)\right),
\end{equation}
and its generating function satisfies
\begin{equation}\label{Sw}
    S_{\mathcal W}(x)=x\left(\log\left(\frac1{1-S_{\mathcal W}(x)}\right)+1\right).
\end{equation}
Expanding  both sides of the equation and identifying the coefficients, we get a system, the resolution of which allows us to obtain the first cases of the enumeration:
\[
1, 2, 9, 68, 730, 10164, 173838, 3524688, 82627200,\dots\quad \href{http://oeis.org/A038037 }{A000169}\quad \cite{Sloane}
\]
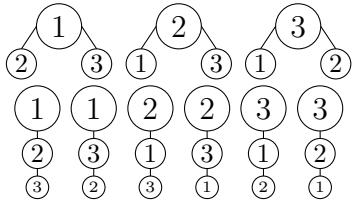
\begin{figure}[ht]
\begin{center}
\begin{tikzpicture}
\draw (0,0) circle (0.3cm);
\draw (-0.5,-0.5) circle (0.2cm);
\draw (0.5,-0.5) circle (0.2cm);
\draw (-0.3,0) -- (-0.5,-0.3);
\draw (0.3,0) -- (0.5,-0.3);
\node at (0,0) {$1$}; 
\node at (-0.5,-0.5) {\footnotesize$2$}; 
\node at (0.5,-0.5) {\footnotesize$3$};
\end{tikzpicture}
\begin{tikzpicture}
\draw (0,0) circle (0.3cm);
\draw (-0.5,-0.5) circle (0.2cm);
\draw (0.5,-0.5) circle (0.2cm);
\draw (-0.3,0) -- (-0.5,-0.3);
\draw (0.3,0) -- (0.5,-0.3);
\node at (0,0) {$2$}; 
\node at (-0.5,-0.5) {\footnotesize$1$}; 
\node at (0.5,-0.5) {\footnotesize$3$};
\end{tikzpicture}
\begin{tikzpicture}
\draw (0,0) circle (0.3cm);
\draw (-0.5,-0.5) circle (0.2cm);
\draw (0.5,-0.5) circle (0.2cm);
\draw (-0.3,0) -- (-0.5,-0.3);
\draw (0.3,0) -- (0.5,-0.3);
\node at (0,0) {$3$}; 
\node at (-0.5,-0.5) {\footnotesize$1$}; 
\node at (0.5,-0.5) {\footnotesize$2$};
\end{tikzpicture}\ \\
\begin{tikzpicture}
\draw (0,0) circle (0.3cm);
\draw (0,-0.6) circle (0.2cm);
\draw (0,-1.05) circle (0.15cm);
\draw (0,-0.3) -- (0,-0.4);
\draw (0,-0.8) -- (0,-0.9);
\node at (0,0) {$1$}; 
\node at (0,-0.6) {\footnotesize$2$}; 
\node at (0,-1.05) {\tiny$3$};
\end{tikzpicture}
\begin{tikzpicture}
\draw (0,0) circle (0.3cm);
\draw (0,-0.6) circle (0.2cm);
\draw (0,-1.05) circle (0.15cm);
\draw (0,-0.3) -- (0,-0.4);
\draw (0,-0.8) -- (0,-0.9);
\node at (0,0) {$1$}; 
\node at (0,-0.6) {\footnotesize$3$}; 
\node at (0,-1.05) {\tiny$2$};
\end{tikzpicture}
\begin{tikzpicture}
\draw (0,0) circle (0.3cm);
\draw (0,-0.6) circle (0.2cm);
\draw (0,-1.05) circle (0.15cm);
\draw (0,-0.3) -- (0,-0.4);
\draw (0,-0.8) -- (0,-0.9);
\node at (0,0) {$2$}; 
\node at (0,-0.6) {\footnotesize$1$}; 
\node at (0,-1.05) {\tiny$3$};
\end{tikzpicture}
\begin{tikzpicture}
\draw (0,0) circle (0.3cm);
\draw (0,-0.6) circle (0.2cm);
\draw (0,-1.05) circle (0.15cm);
\draw (0,-0.3) -- (0,-0.4);
\draw (0,-0.8) -- (0,-0.9);
\node at (0,0) {$2$}; 
\node at (0,-0.6) {\footnotesize$3$}; 
\node at (0,-1.05) {\tiny$1$};
\end{tikzpicture}
\begin{tikzpicture}
\draw (0,0) circle (0.3cm);
\draw (0,-0.6) circle (0.2cm);
\draw (0,-1.05) circle (0.15cm);
\draw (0,-0.3) -- (0,-0.4);
\draw (0,-0.8) -- (0,-0.9);
\node at (0,0) {$3$}; 
\node at (0,-0.6) {\footnotesize$1$}; 
\node at (0,-1.05) {\tiny$2$};
\end{tikzpicture}
\begin{tikzpicture}
\draw (0,0) circle (0.3cm);
\draw (0,-0.6) circle (0.2cm);
\draw (0,-1.05) circle (0.15cm);
\draw (0,-0.3) -- (0,-0.4);
\draw (0,-0.8) -- (0,-0.9);
\node at (0,0) {$3$}; 
\node at (0,-0.6) {\footnotesize$2$}; 
\node at (0,-1.05) {\tiny$1$};
\end{tikzpicture}
\caption{The nine windmills of degree $3$.\label{windmills}}

\end{center}
\end{figure}

This is another example of a tree-like structure studied among others in \cite{BLL}. Notice that no closed form for the generating function is known but there exists a formula for the coefficients as a combination of Stirling numbers of first kind,
\begin{equation}\label{Wn}
    W_n=\sum_{i=0}^ni!\binom nis_{n-1,i},
\end{equation}
where $s_{n,i}$ denotes the (unsigned) Stirling number of first kind that counts the number of permutations of $n$ objects with exactly $i$ cycles. Indeed, from (\ref{Sw}), $S_{\mathcal W}(x)$ is the inverse of $g(x)={x\over 1-\log(1-x)}$ for the composition. The Lagrange inversion theorem \cite{Lagrange} is a classical combinatorial tools allowing us to compute the Taylor expansion of inverse function. In our case, the direct application of the Lagrange inversion Theorem implies that \begin{equation}\label{lagrange}W_n=\left.\left(d\over dx\right)^{n-1}\left(x\over g(x)\right)^n\right|_{x=0}=\left.\left(d\over dx\right)^{n-1}\left(1+\log\left(1\over1-x\right)\right)^n\right|_{x=0}\end{equation}
is the coefficient of $x^{n-1}$ in $\left(1-\log\left(1\over1-x\right)\right)^n$ multiplied by $(n-1)!$.  Knowing that the exponential generating function of Stirling's numbers $s_{i,k}$  ($k$ fixed) is $\frac1k!\log\left(\frac1{1-x}\right)^k=\sum_{i}s_{i,k}{x^i\over i!}$, an easy computation allows us to deduce (\ref{Wn}) from (\ref{lagrange}).
\\ \\
 The last example we consider is the one where no more order constraints are imposed on the sub-trees of the same node. In this context, a tree is a root with a (possibly empty) set of trees. The trees of this kind can be drawn as  nesting of disjointed discs with numbered surfaces (see some examples in Figure \ref{nested}). Notice that, in these examples nested disk configurations, of degree $3$ are as numerous as the windmills of degree $3$ (see Figure \ref{windmills}). 
 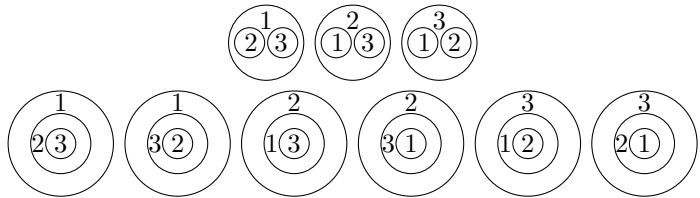
\begin{figure}[ht]
 \begin{center}
     \begin{tikzpicture}
     \draw(0,0) circle (0.5cm);
     \draw(-0.22,0) circle (0.2cm);
     \draw(0.22,0) circle (0.2cm);
     \node at (0,0.3) {\footnotesize$1$};
     \node at (-0.2,0) {\footnotesize$2$};
     \node at (0.2,0) {\footnotesize$3$};
     \end{tikzpicture}
     \begin{tikzpicture}
     \draw(0,0) circle (0.5cm);
     \draw(-0.22,0) circle (0.2cm);
     \draw(0.22,0) circle (0.2cm);
     \node at (0,0.3) {\footnotesize$2$};
     \node at (-0.2,0) {\footnotesize$1$};
     \node at (0.2,0) {\footnotesize$3$};
     \end{tikzpicture}
    \begin{tikzpicture}
     \draw(0,0) circle (0.5cm);
     \draw(-0.22,0) circle (0.2cm);
     \draw(0.22,0) circle (0.2cm);
     \node at (0,0.3) {\footnotesize$3$};
     \node at (-0.2,0) {\footnotesize$1$};
     \node at (0.2,0) {\footnotesize$2$};
     \end{tikzpicture}\\
        \begin{tikzpicture}
     \draw(0,0) circle (0.7cm);
     \draw(0,0) circle (0.2cm);
     \draw(0,0) circle (0.4cm);
     \node at (0,0.55) {\footnotesize$1$};
     \node at (-0.3,0) {\footnotesize$2$};
     \node at (0,0) {\footnotesize$3$};
     \end{tikzpicture}
      \begin{tikzpicture}
     \draw(0,0) circle (0.7cm);
     \draw(0,0) circle (0.2cm);
     \draw(0,0) circle (0.4cm);
     \node at (0,0.55) {\footnotesize$1$};
     \node at (-0.3,0) {\footnotesize$3$};
     \node at (0,0) {\footnotesize$2$};
     \end{tikzpicture}
      \begin{tikzpicture}
     \draw(0,0) circle (0.7cm);
     \draw(0,0) circle (0.2cm);
     \draw(0,0) circle (0.4cm);
     \node at (0,0.55) {\footnotesize$2$};
     \node at (-0.3,0) {\footnotesize$1$};
     \node at (0,0) {\footnotesize$3$};
     \end{tikzpicture}
      \begin{tikzpicture}
     \draw(0,0) circle (0.7cm);
     \draw(0,0) circle (0.2cm);
     \draw(0,0) circle (0.4cm);
     \node at (0,0.55) {\footnotesize$2$};
     \node at (-0.3,0) {\footnotesize$3$};
     \node at (0,0) {\footnotesize$1$};
     \end{tikzpicture}
      \begin{tikzpicture}
     \draw(0,0) circle (0.7cm);
     \draw(0,0) circle (0.2cm);
     \draw(0,0) circle (0.4cm);
     \node at (0,0.55) {\footnotesize$3$};
     \node at (-0.3,0) {\footnotesize$1$};
     \node at (0,0) {\footnotesize$2$};
     \end{tikzpicture}
      \begin{tikzpicture}
     \draw(0,0) circle (0.7cm);
     \draw(0,0) circle (0.2cm);
     \draw(0,0) circle (0.4cm);
     \node at (0,0.55) {\footnotesize$3$};
     \node at (-0.3,0) {\footnotesize$2$};
     \node at (0,0) {\footnotesize$1$};
     \end{tikzpicture}
     \end{center}
     \caption{The nine configurations of nested discs of degree $3$. \label{nested}}
 \end{figure}
 Of course, this is not always the case and generally, there are fewer  nested discs configurations  than windmills. For instance, there are two windmills of degree $4$ the roots of which is labeled by $1$ with three sub-windmills of degree $1$ while there is only one nested discs configuration the root of which is labeled by $1$ containing three discs (see Figure \ref{wvsn}).
 \begin{figure}[ht]
 \begin{center}
     \begin{tikzpicture}
     \draw(0,0) circle (0.4cm);
     \draw(-0.5,-0.7) circle (0.2cm);
     \draw(0,-0.7) circle (0.2cm);
     \draw(0.5,-0.7) circle (0.2cm);
     \draw (-0.4,0)--(-0.5,-0.5);
     \draw (0.4,0)--(0.5,-0.5);
     \draw (0,-0.4)--(0,-0.5);
     \node at(0,0) {$1$};
     \node at(-0.5,-0.7) {\footnotesize$2$};
     \node at(0,-0.7) {\footnotesize$3$};
     \node at(0.5,-0.7) {\footnotesize$4$};
     \end{tikzpicture}
     \begin{tikzpicture}
     \draw(0,0) circle (0.4cm);
     \draw(-0.5,-0.7) circle (0.2cm);
     \draw(0,-0.7) circle (0.2cm);
     \draw(0.5,-0.7) circle (0.2cm);
      \draw (-0.4,0)--(-0.5,-0.5);
     \draw (0.4,0)--(0.5,-0.5);
     \draw (0,-0.4)--(0,-0.5);
     \node at(0,0) {$1$};
     \node at(-0.5,-0.7) {\footnotesize$2$};
     \node at(0,-0.7) {\footnotesize$4$};
     \node at(0.5,-0.7) {\footnotesize$3$};
     \end{tikzpicture}
         \begin{tikzpicture}
     \draw(0,0) circle (0.9cm);
     \draw(-0.5,-0.4) circle (0.2cm);
     \draw(0,0.6) circle (0.2cm);
     \draw(0.5,-0.4) circle (0.2cm);
     \node at(0,0) {$1$};
     \node at(-0.5,-0.4) {\footnotesize$2$};
     \node at(0,0.6) {\footnotesize$3$};
     \node at(0.5,-0.4) {\footnotesize$4$};
     \end{tikzpicture}
     \caption{Two windmills and one nested discs configuration. \label{wvsn}}
 \end{center}
 \end{figure}
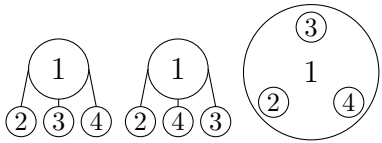
 \\
 The combinatorial class satisfies the isomorphism
 \[
 \mathcal{N}pt=\bullet\boxtimes\Set(\mathcal Npt),
 \]
 and its generating series satisfies the functional equation
 \[
 S_{\mathcal Npt}(x)=xe^{S_{\mathcal Npt}(x)}.
 \]
 Solving these equation, one finds 
 \[
  S_{\mathcal Npt}(x)=-W(-x),
 \]
 where $W(x)$ denotes the Lambert $W$ function that is the principal branch of the functional inverse of $x\rightarrow xe^x$ \cite{CGHJK}. The Taylor expansion of $W(x)$ is obtained by applying Lagrange inversion Theorem and implies $Npt_n=n^{n-1}$.
 The sequence of the $Npt_n$'s can also be found in \cite{Sloane}:
 \[
 1, 2, 9, 64, 625, 7776, 117649, 2097152,\dots\quad\href{https://oeis.org/A000169}{A000169}.
 \]

\end{document}